\documentclass[12pt]{article}

\usepackage{amsmath}
\usepackage{amssymb}
\usepackage{amscd}
\usepackage{amsthm}
\theoremstyle{plain}
\newtheorem{Thm}{Theorem}[section]

\newtheorem{Prop}[Thm]{Proposition}

\theoremstyle{definition}
\newtheorem{Defn}[Thm]{Definition}
\theoremstyle{remark}
\newtheorem{Rem}[Thm]{Remark}
\numberwithin{equation}{section}

\def\Rat{\mathop{\mathrm{Rat}}}

\def\Hilb{\mathop{\mathrm{Hilb}}}

\begin{document}
\title{A note on the projective varieties \\ of almost general type}
\author{ Shigetaka FUKUDA}

\date{\empty}

\maketitle \thispagestyle{empty}
\pagestyle{myheadings}
\markboth{Shigetaka FUKUDA}{The projective varieties of almost general type}
\begin{abstract}
A $\mathbf{Q}$-Cartier divisor $D$ on a projective variety $M$ is {\it almost nup}, if $(D , C) > 0$ for every very general curve $C$ on $M$.
An algebraic variety $X$ is of {\it almost general type}, if there exists a projective variety $M$ with only terminal singularities such that the canonical divisor $K_M$ is almost nup and such that $M$ is birationally equivalent to $X$.
We prove that a complex algebraic variety is of almost general type if and only if it is neither uniruled nor covered by any family of varieties being birationally equivalent to minimal varieties with numerically trivial canonical divisors, under the minimal model conjecture.
Furthermore we prove that, for a projective variety $X$ with only terminal singularities, $X$ is of almost general type if and only if the canonical divisor $K_X$ is almost nup, under the minimal model conjecture.
\end{abstract}

\section{Introduction}

Throughout the paper every variety is defined over the field of complex numbers $\mathbf{C}$.
We follow the notation and terminology in \cite{Utah}.

\begin{Defn}
A $\mathbf{Q}$-Cartier divisor $D$ on a projective variety $X$ is {\it almost numerically positive} ({\it almost nup}, for short), if there exists a union $F$ of at most countably many prime divisors on $X$ such that $(D , C) > 0$ for every curve $C \nsubseteq F$ (i.e.\ if $(D , C) > 0$ for every very general curve $C$).
We say that $D$ is {\it quasi-numerically positive} ({\it quasi-nup}, for short), if $D$ is nef and almost nup.
\end{Defn}

\begin{Defn}
An algebraic variety $X$ is {\it of almost general type}, if there exists a projective variety $M$ with only terminal singularities such that the canonical divisor $K_M$ is almost nup and that $M$ is birationally equivalent to $X$.
\end{Defn}

Obviously, an algebraic variety $X$ is of almost general type, if it is of general type (i.e.\ if the geometric Kodaira dimension $\kappa_{geom} (X)$ equals the dimension $\dim X$).
However, the converse statement contains the conjecture that every quasi-nup canonical divisor on a minimal variety should be semi-ample, which is the essential part of the famous abundance conjecture (cf.\ Ambro \cite{Am}).

\begin{Defn}
The property $P(x)$ holds {\it for every very general point} $x$ on an algebraic variety $X$, if $P(x)$ holds for every point $x \in X \setminus E$ where $E$ is some fixed union of at most countably many prime divisors on $X$.
\end{Defn}

\begin{Defn}
For an ambient space $X$, a subset $L$ is {\it covered} by subsets $D_i$ ($i \in I$) if $L \subset \bigcup_{i \in I} D_i$.
\end{Defn}

\begin{Defn}
For an algebraic variety $X$, we denote by $\Rat X$ the rational function field of $X$.
\end{Defn}

In this article, we look into the complex algebraic varieties of almost general type and show that the class of varieties of this type coincides with the class of varieties that are neither uniruled nor covered by any family of varieties being birationally equivalent to minimal varieties with numerically trivial canonical divisors, under the minimal model conjecture (see Theorem \ref{Thm:MT}).
In the category of projective varieties with only terminal singularities, this class is proved to be the class of varieties with almost nup canonical divisors, under the minimal model conjecture (see Theorem \ref{Thm:Iff}).
In contrast with the method in the previous paper \cite{Fuk}, we proceed not depending on the (log) abundance conjecture, which claims the semi-ampleness of nef (log) canonical divisors, and the deformation invariance theorem of Kodaira dimension due to Tsuji \cite{Ts02} and Siu \cite{Si}.
The key is Tsuji's theory \cite{Tsuji} of numerically trivial fibrations.

We briefly explain the background.
The minimal model conjecture states that every projective variety with only $\mathbf{Q}$-factorial terminal singularities would reach a variety with nef canonical divisor (i.e.\ a minimal variety) or with some Mori fiber space structure, after a finite sequence of divisorial contractions and flips.
We emphasize that this conjecture does not contain the abundance conjecture that the (nef) canonical divisor of every minimal variety should be semi-ample.
Compared with the mechanism of the minimal model theory, the known proofs of the abundance theorems for surfaces and threefolds are case-by-case and very complicated.
So it is meaningful to get geometric results, not depending on the abundance conjecture, but only on the minimal model conjecture.

\begin{Rem}
We state the current status of the minimal model conjecture for the convenience of the reader.
The general theory reduces this conjecture to the existence and termination of flips.
The existence is now a theorem due to Birkar, Cascini, Hacon and McKernan \cite{BCHM} in all dimensions.
The termination became a theorem (Shokurov \cite{Sh85} for threefolds and Kawamata-Matsuda-Matsuki \cite{KMM} for fourfolds) in dimension $\leq 4$.
\end{Rem}

\section{Stability under deformation}

In this section, we show that the property ``being of almost general type'' is stable under deformation.

\begin{Prop}\label{Prop:Stab1}
Let $f:M \to S$ be a surjective morphism between smooth projective varieties such that the field extension $\Rat M / \Rat S$ is algebraically closed.
Then one of the following holds:

(i) The canonical divisor $K_F$ is not almost nup for every general fiber $F$ of $f$.

(ii) The canonical divisor $K_F$ is almost nup for every very general fiber $F$ of $f$.
\end{Prop}

\begin{proof}
Let $\mathcal{H} \subset M \times \Hilb (M)$ be the universal family parametrized by the Hilbert scheme $\Hilb (M)$.
Consider the projection morphisms $p_1 : \mathcal{H} \to M$ and $p_2 : \mathcal{H} \to \Hilb (M)$.

Let $H$ be a hyperplane section of $S$.
The morphism $p_2$ is flat.
Thus, for every curve $C$ on $M$, there exists some irreducible component 
$V$ of $\mathcal{H}$ such that $V \supseteq p_2^{-1} ( [C] )$ and that $\dim 
V - \dim p_2 (V) =1$, where $[C]$ is the point representing the subscheme 
$C$ of $M$.
Now define the set $I:= \{ V \vert V$ is an irreducible component of $\mathcal{H}$, such that $V \supseteq p_2^{-1} ( [C] )$ for some curve $C$ with $(K_M, C) \leq 0$ and $(f^* H, C) = 0$ and  that $\dim V - \dim p_2 (V) =1 \}$.
Note that $(f^* H, C) = 0$ if and only if $f(C)$ is a point.

We divide the situation into two cases, according to the bigness of the set $\bigcup \{p_1 (V) \vert V \in I \}$ in $M$.

First we treat the case where $p_1 (V) = M$ for some $V \in I$.

We apply the normalization, the Stein factorization and the flattening to the morphism $p_2 \vert_V: V \to p_2 (V)$ and obtain the following commutative diagram among projective varieties
 \begin{equation}
  \begin{CD}
   W @>\text{$\mu$}>> T \\
   @V\text{$\rho$}VV @VV\text{$\nu$}V \\
   V @>>\text{$p_2 \vert_V$}> p_2(V)
  \end{CD}
 \end{equation}
with the properties:

(1) the morphism $\rho$ is birational

(2) the morphism $\nu$ is generically finite

(3) the extension $\Rat W / \Rat T$ is algebraically closed

(4) the morphism $\mu$ is flat

(5) the variety $T$ is smooth.

For every curve $C' \subseteq (\nu \mu)^{-1} ([C]) = \rho^{-1} (p_2^{-1} ([C]))$, the relation $((p_1 \rho)^* K_M, C') \leq 0$ and $((f p_1 \rho)^* H, C') = 0$ holds.
Thus for every fiber $F$ of $\mu$, we have the relation $((p_1 \rho)^* K_M, F) \leq 0$ and $((f p_1 \rho)^* H, F) = 0$, because all fibers of the flat morphism $\mu$ represent the same homology class in $H_2 (W, \mathbf{Z})$.

Let $W'$ be the image of the morphism $(p_1 \rho , \mu): W \to M \times T$.
We note that $W'$ ($\subseteq M \times T$) is birationally equivalent to $W$, because the image of the morphism $(id., \nu)(p_1 \rho , \mu): W \to M \times T \to M \times p_2 (V)$ is $V$.
Consider the projection morphisms $p_1' : W' \to M$ and $p_2' : W' \to T$.
Then $p_1'(W') = M$.
\begin{equation}
 \begin{CD}
W \\
@VVV \\
W' @>\text{embedding}>> M \times T @>\text{2nd proj.}>> T \\
@.   @V\text{1st proj.}VV   \\
@.   M  \\
@.   @V\text{$f$}VV  \\
@.   S
 \end{CD}
\end{equation}

We have some open dense subset $U$ of $T$ such that every fiber $G$ of $p_2'$ over $U$ is irreducible.
By applying the projection formula to the birational morphism $W \to W'$, the relation that $((p_1 ')^* K_M, G) \leq 0$ and $((f p_1 ')^* H, G) = 0$ follows.
Note that $G$ can be identified with a curve on $M$ that is $f$-exceptional and whose intersection number with $K_M$ is not positive.
The property of constructible sets implies that $p_1 '((p_2 ')^{-1} (U))$ and $f p_1 '((p_2 ')^{-1} (U))$ contain open dense subsets of $M$ and $S$, respectively.

Consequently the statement (i) holds.

Next we treat the case where $p_1 (V) \subsetneq M$ for every $V \in I$.
The countability of components of the Hilbert scheme $\Hilb (M)$ implies that the statement (ii) holds.
\end{proof}

We cite Tsuji's existence theorem of numerically trivial fibrations.

\begin{Prop}[Tsuji \cite{Tsuji}, cf.\ \cite{BCEKPRSW}]\label{Prop:Tsuji}
Let $X$ be a normal projective variety and $L$ a nef divisor on $X$.
Then there exist projective varieties $Y$ and $Z$ and morphisms $\mu: Y \to X$ and $\nu: Y \to Z$ with the following properties:

(i) the morphism $\mu$ is birational

(ii) the morphism $\nu$ is surjective

(iii) the extension $\Rat Y / \Rat Z$ is algebraically closed

(iv) for some two open subvarieties $U$ and $V$ of $X$ and $Z$ respectively, $\mu^{-1} (U) = \nu^{-1} (V)$ and $\mu \vert_{\mu^{-1} (U)}$ is isomorphic

(v) the divisor $\mu^* L \vert_F$ is numerically trivial for every very general fiber $F$ of $\nu$

(vi) for every very general point $x \in X$ there does not exist a closed subvariety $S \ni x$ such that $L \vert_S$ is numerically trivial and that $\dim S > \dim Y - \dim Z$.
\end{Prop}

The following means that the property for a canonical divisor to be almost nup is stable under birational transformation.

\begin{Thm}\label{Thm:Iff}
Assume that the minimal model conjecture holds in dimension $n$.
Let $X$ be a projective variety with only terminal singularities of dimension $n$.
Then the canonical divisor $K_X$ is almost nup if and only if $X$ is of almost general type.
\end{Thm}

\begin{proof}
The ``only if'' part.
Trivial from the definition of ``being of almost general type''.

The ``if'' part.
From Miyaoka-Mori (\cite{MiMo}), $X$ is not uniruled.
Thus a crepant $\mathbf{Q}$-factorialization of $X$ has a minimal model $Z$.
Here we have a common resolution $\mu: Y \to X$ and $\nu: Y \to Z$ such that $\mu^{*} K_X = \nu^{*} K_Z + E$ where $E$ is a $\nu$-exceptional effective $\mathbf{Q}$-divisor.
Suppose that $K_Z$ is not quasi-nup.
Then for some open subvariety $U$ of $Z$ there exists a proper dominating morphism $\rho$ from $U$ to a lower-dimensional variety $V$ such that the extension $\Rat U / \Rat V$ is algebraically closed and that $K_Z \vert_F$ is numerically trivial for every very general fiber $F$ of $\rho$, from Proposition \ref{Prop:Tsuji} (Tsuji).
For every desingularization $\alpha:M \to Z$, the divisor $(K_M - \alpha^* K_Z) \vert_{\alpha^{-1} (F)}$ is $\alpha \vert_{\alpha^{-1} (F)}$-exceptional.
Thus $(K_M - \alpha^* K_Z) \alpha^* (H)^{\dim F -1} \vert_{\alpha^{-1} (F)} = 0 $ where $H$ is a hyperplane section of $Z$.
So $K_M \alpha^* (H)^{\dim F -1} \vert_{\alpha^{-1} (F)} = 0 $.
Therefore $K_M$ is not almost nup and thus $X$ is not of almost general type.
This is a contradiction!
Consequently $K_Z$ is quasi-nup and thus $\mu^* K_X$ is almost nup.
As a result, $K_X$ is almost nup.
\end{proof}

Here we have the main result of this section.

\begin{Thm}\label{Thm:Stab2}
Assume that the minimal model conjecture holds in dimension $n$.
Let $f:M \to S$ be a surjective morphism between (possibly singular) projective varieties with relative dimension $n$ such that the extension $\Rat M / \Rat S$ is algebraically closed.
Then one of the following holds:

(i) Every general fiber $F$ of $f$ is not of almost general type.

(ii) Every very general fiber $F$ of $f$ is of almost general type.
\end{Thm}

\begin{proof}
Proposition \ref{Prop:Stab1} and Theorem \ref{Thm:Iff} imply the assertion.
\end{proof}

\section{Some kind of hyperbolicity}

In this section, we show that the varieties of almost general type are characterized by some kind of hyperbolicity.

\begin{Thm}\label{Thm:Hyp1}
Assume that the minimal model conjecture holds in dimension $< n$.
Let $X$ be a projective variety with only terminal singularities of dimension $n$.
If $K_X$ is almost nup, then the locus $\bigcup \{ D ; \thickspace D$ is a closed subvariety ($\subsetneqq X$) not of almost general type $\} $ is covered by at most countably many prime divisors on $X$.
\end{Thm}

\begin{proof}
The proof proceeds along the same line as in the paper \cite{Fuk}, by using Theorem \ref{Thm:Stab2}.
But, for the readers' convenience, we do not make the presentation rough.

Assuming that $K_X$ is almost nup and that however the locus $\bigcup \{ D ; \thickspace D$ is a closed subvariety ($\subsetneqq X$) not of almost general type $\}$ cannot be covered by at most countably many prime divisors on $X$, we derive a contradiction.

Let $\mathcal{H} \subset X \times \Hilb (X)$ be the universal family parametrized by the Hilbert scheme $\Hilb (X)$.
By the countability of the components of $\Hilb (X)$, we have an irreducible component $V$ of $\mathcal{H}$ with surjective projection morphisms $f:V \to X$ and $g:V \to T ( \subset \Hilb (X))$ from $V$ to projective varieties $X$ and $T$ respectively, such that $f(g^{-1} (t)) \subsetneqq X$ for every $t \in T$ and that the locus $\bigcup \{ D ; \thickspace D$ is a closed subvariety ($\subsetneqq X$) not of almost general type and $D = f(g^{-1} (t))$ for some $t \in T \}$ cannot be covered by at most countably many prime divisors on $X$.
Let $\nu : V_{norm} \to V$ be the normalization.
We consider the Stein factorization of $g \nu$ into the finite morphism $g_1 : S \to T$ from a projective normal variety $S$ and the morphism $g_2 : V_{norm} \to S$ with an algebraically closed extension $\Rat V / \Rat S$.
Put $V^* :=$ [the image of the morphism $(f \nu, g_2) : V_{norm} \to X \times S$].
\begin{equation}
 \begin{CD}
V_{norm} @>>> V^* @>\text{embedding}>> X \times S @>>> S    \\
@.   @VVV @VVV @VV\text{$g_1$}V \\
@.   V   @>\text{embedding}>> X \times T @>>> T    \\
@.   @.   @VVV \\
@.   @.   X
 \end{CD}
\end{equation}
Note that every fiber of the morphism $g:V \to T$ consists of a finite number of fibers of the projection morphism from $V^*$ to $S$.
Thus we may replace $(V,T)$ by $(V^* , S)$ and assume that the extension $\Rat V / \Rat T$ is algebraically closed.

According to Theorem \ref{Thm:Stab2} (i) and (ii), we divide the situation into two cases.

First consider the case where $g^{-1} (t)$ is of almost general type for very general $t \in T$.
Then there exists a subvariety $T_0 \subsetneqq T$ such that the locus $\bigcup \{ D ; \thickspace D$ is a closed subvariety ($\subsetneqq X$) not of almost general type and $D = f(g^{-1} (t))$ for some $t \in T_0 \}$ cannot be covered by at most countably many prime divisors on $X$.
Thus we can replace $(V, T)$ by $(V_1 , T_1)$, where $V_1$ and $T_1$ are projective varieties such that $V_1$ is some suitable irreducible component of $g^{-1} (T_0)$ and $T_1 = g(V_1)$ .
Because $\dim V_1 < \dim V$, by repeating this process of replacement, we can reduce the assertion to the next case.

Now consider the case where $g^{-1} (t)$ is not of almost general type for general $t \in T$.
If $\dim V > \dim X$, then $\dim f^{-1} (x) \geq 1$ for all $x \in X$, thus $\dim g(f^{-1} (x)) \geq 1$ (this means that $g(f^{-1} (x)) \cap H \ne \emptyset$ for every hyperplane section $H$ of $T$).
Therefore, by repeating the process of cutting $T$ by general hyperplanes, we can reduce the assertion to the subcase where $\dim V = \dim X$.
Take birational morphisms $\alpha : X' \to X$ and $\beta : V' \to V$ from non-singular projective varieties with a generically finite morphism $\gamma : V' \to X'$ such that $\alpha \gamma = f \beta$.
\begin{equation}
 \begin{CD}
V @<\text{$\beta$}<< V' \\
@V\text{$f$}VV @VV\text{$\gamma$}V \\
X @<<\text{$\alpha$}< X'
 \end{CD}
\end{equation}
Because $K_X$ is almost nup and $K_{X'} - \alpha^* K_X$ is effective, $K_{X'}$ is almost nup.
Here we have $K_{V'} = \gamma^{*} K_{X'} + R_{\gamma}$, where $R_{\gamma}$ is the ramification divisor (which is effective) for $\gamma$.
Thus $K_{V'}$ becomes almost nup.
For a very general fiber $F' = (g \beta )^{-1} (t) $ for $g \beta$ (i.e.\ the point $t$ does not belong to some fixed union of at most countably many prime divisors on $T$), also $K_{F'}$ is almost nup.
Thus every very general fiber $F$ of $g$ is of almost general type.
This is a contradiction!
\end{proof} 

As a corollary we have

\begin{Thm}
Let $X$ be a projective variety with only terminal singularities of dimension $\leq 5$.
If $K_X$ is almost nup, then the locus $\bigcup \{ D ; \thickspace D$ is a closed subvariety ($\subsetneqq X$) not of almost general type $\} $ is covered by at most countably many prime divisors on $X$.
\end{Thm}

\begin{proof}
Theorem \ref{Thm:Hyp1} implies the assertion by virtue of the existence theorem (Shokurov \cite{Sh}, cf.\ Hacon and McKernan \cite{HM}) and the termination theorem (Kawamata-Matsuda-Matsuki \cite{KMM}) of four-dimensional flips. 
\end{proof}

We cite the following proposition and proof for the convenience of the reader.

\begin{Prop}[\cite{Fuk}]\label{Prop:NonPar}
Let $X$ be a projective variety with only terminal singularities.
If $K_X$ is almost nup, then $X$ is not birationally equivalent to any minimal variety with numerically trivial canonical divisor.
\end{Prop}

\begin{proof}[Proof (\cite{Fuk})]
We give a sketchy proof.
Let $Z$ be a minimal variety with numerically trivial canonical divisor $K_Z$.
For any resolution $\nu:Y \to Z$, the divisor $K_Y-\nu^*K_Z$ is effective and $\nu$-exceptional.
Thus $K_Y (\nu^* H)^{\dim Y -1} =0$ for any ample divisor $H$ on $Z$.
Consequently $K_Y$ is not almost nup.
\end{proof}

Now we consider the converse statement of Theorem \ref{Thm:Hyp1} and Proposition \ref{Prop:NonPar}.

\begin{Thm}\label{Thm:Conv}
Assume that the minimal model conjecture holds in dimension $n$.
Let $X$ be a projective variety with only terminal singularities of dimension $n$.
If $X$ is neither uniruled nor birationally equivalent to any variety fibered by minimal varieties with numerically trivial canonical divisors, then $K_X$ is almost nup.
\end{Thm}

\begin{proof}
A crepant $\mathbf{Q}$-factorialization of $X$ has a minimal model $Z$.
Suppose that $K_X$ is not almost nup.
From the inequality $K_X \geq K_Z$, between bi-divisors (\cite{Sh}), $K_Z$ is not quasi-nup.
Thus from Proposition \ref{Prop:Tsuji} (Tsuji), $Z$ is birationally equivalent to some variety fibered by minimal varieties with numerically trivial canonical divisors.
This is a contradiction!
Consequently $K_X$ is almost nup.
\end{proof}

At last we establish the following

\begin{Thm}[Main Theorem]\label{Thm:MT}
Assume that the minimal model conjecture holds in dimension $\leq n$.
Let $X$ be a projective variety with only terminal singularities of dimension $n$.
Then the five conditions below are equivalent to each other:

(1) $X$ is of almost general type.

(2) $K_X$ is almost nup.

(3) The locus $\bigcup \{ D ; \thickspace D$ is a closed subvariety ($\subsetneqq X$) not of almost general type $\} $ is covered by at most countably many prime divisors on $X$ and the variety $X$ is not birationally equivalent to any minimal variety with numerically trivial canonical divisor and is not a rational curve.

(4) $X$ is not uniruled and can not be covered by any family of varieties being birationally equivalent to minimal varieties with numerically trivial canonical divisors.

(5) $X$ is neither uniruled nor birationally equivalent to any variety fibered by minimal varieties with numerically trivial canonical divisors.
\end{Thm}

\begin{proof}
(1) and (2) are equivalent (Theorem \ref{Thm:Iff}).
(2) implies (3) by Theorem \ref{Thm:Hyp1} and Proposition \ref{Prop:NonPar}.
(3) implies (4), because every minimal variety with numerically trivial canonical divisor is not of almost general type from Theorem \ref{Thm:Iff}.
Obviously (4) implies (5).
(5) implies (2) by Theorem \ref{Thm:Conv}.
\end{proof}

\begin{Rem}[cf.\ \cite{Fuk}]
The assumption of Main Theorem \ref{Thm:MT} is satisfied for $n \leq 4$, because of the existence (Shokurov \cite{Sh}, cf.\ Hacon and McKernan \cite{HM}) and the termination (Kawamata-Matsuda-Matsuki \cite{KMM}) of four-dimensional flips. 
\end{Rem}

\section*{Appendix}

Here we note the theory of pseudo-effective divisors due to Boucksom-Demailly-Paun-Peternell (\cite{BDPP}).
Let $L$ be a Cartier divisor on a smooth projective variety $X$.
The divisor $L$ is pseudo-effective if and only if $(L,C) \geq 0$ for every curve $C$ through a very general point $x$ on $X$ (\cite[Theorem 0.2]{BDPP}).
For the pseudo-effective divisor $L$, we set the number $p(L):= \dim X$ (the notation at \cite[Definition 8.3]{BDPP}) if and only if $L$ is almost nup.

If $X$ is of dimension $4$, the canonical divisor $K_X$ is pseudo-effective and $p(K_X) = 4$, then every proper subvariety $S \subset X$ through a very general point $x$ on $X$ is of general type.
(The paper \cite{Fuk} overlaps this result \cite[Proposition 9.12]{BDPP}. The presentation of the proof in the latter paper is casual.)

\bigskip
Faculty of Education, Gifu Shotoku Gakuen University

Yanaizu-cho, Gifu City, Gifu 501-6194, Japan

fukuda@ha.shotoku.ac.jp

\end{document}